\documentclass[11pt]{article}

\usepackage{latexsym}
\usepackage{graphicx}
\usepackage{framed}
\usepackage{threeparttable}
\usepackage{amssymb}
\usepackage{amsmath}
\usepackage{multirow}
\usepackage{hyperref}

\begin{document}

\newcommand{\noi}{\noindent}
\newcommand{\nn}{\nonumber}
\newcommand{\bd}{\begin{displaymath}}
\newcommand{\ed}{\end{displaymath}}
\newcommand{\bp}{\underline{\bf Proof}:\ }
\newcommand{\ep}{{\hfill $\Box$}\\ }
\newtheorem{1}{LEMMA}[section]
\newtheorem{2}{THEOREM}[section]
\newtheorem{3}{COROLLARY}[section]
\newtheorem{4}{PROPOSITION}[section]
\newtheorem{5}{REMARK}[section]
\newtheorem{20}{OBSERVATION}[section]
\newtheorem{10}{DEFINITION}[section]
\newtheorem{30}{RESULTS}[section]
\newtheorem{40}{CLAIM}[section]
\newtheorem{50}{ASSUMPTION}[section]
\newtheorem{60}{EXAMPLE}[section]
\newtheorem{70}{ALGORITHM}[section]
\newtheorem{80}{PROBLEM}
\newcommand{\be}{\begin{equation}}
\newcommand{\ee}{\end{equation}}
\newcommand{\ba}{\begin{array}}
\newcommand{\ea}{\end{array}}
\newcommand{\bea}{\begin{eqnarray}}
\newcommand{\eea}{\end{eqnarray}}
\newcommand{\bqn}{\begin{eqnarray*}}
\newcommand{\eqn}{\end{eqnarray*}}


\newcommand{\e} { \ = \ }
\newcommand{\leqs}{ \ \leq \ }
\newcommand{\geqs}{ \ \geq \ }
\def\theequation{\thesection.\arabic{equation}}
\def\bReff#1{{\bRm
(\bRef{#1})}}
\newcommand{\eps}{\varepsilon}
\newcommand{\sgn}{\operatorname{sgn}}
\newcommand{\sign}{\operatorname{sign}}
\newcommand{\Vol}{\operatorname{Vol}}
\newcommand{\Var}{\operatorname{Var}}
\newcommand{\Cov}{\operatorname{Cov}}
\newcommand{\vol}{\operatorname{vol}}
\newcommand{\var}{\operatorname{var}}
\newcommand{\cov}{\operatorname{cov}}
\renewcommand{\Re}{\operatorname{Re}}
\renewcommand{\Im}{\operatorname{Im}}
\newcommand{\bE}{{\mathbb E}}
\newcommand{\bR}{\mathbb{R}}
\newcommand{\bN}{{\mathbb N}}
\newcommand{\bC}{\mathbb{C}}
\newcommand{\bF}{\mathbb{F}}
\newcommand{\bQ}{{\mathbb Q}}
\newcommand{\bZ}{{\mathbb Z}}
\newcommand{\cA}{{\mathcal A}}
\newcommand{\cB}{{\mathcal B}}
\newcommand{\cC}{{\mathcal C}}
\newcommand{\cD}{{\mathcal D}}
\newcommand{\cE}{{\mathcal E}}
\newcommand{\cF}{{\mathcal F}}
\newcommand{\cG}{{\mathcal G}}
\newcommand{\cH}{{\mathcal H}}
\newcommand{\cI}{{\mathcal I}}
\newcommand{\cJ}{{\mathcal J}}
\newcommand{\cK}{{\mathcal K}}
\newcommand{\cL}{{\mathcal L}}
\newcommand{\cM}{{\mathcal M}}
\newcommand{\cN}{{\mathcal N}}
\newcommand{\cO}{{\mathcal O}}
\newcommand{\cP}{{\mathcal P}}
\newcommand{\cQ}{{\mathcal Q}}
\newcommand{\cR}{{\mathcal R}}
\newcommand{\cS}{{\mathcal S}}
\newcommand{\cT}{{\mathcal T}}
\newcommand{\cU}{{\mathcal U}}
\newcommand{\cV}{{\mathcal V}}
\newcommand{\cW}{{\mathcal W}}
\newcommand{\cX}{{\mathcal X}}
\newcommand{\cY}{{\mathcal Y}}
\newcommand{\cZ}{{\mathcal Z}}
\newcommand{\bx}{{\mathbf x}}
\newcommand{\by}{{\mathbf y}}
\newcommand{\bz}{{\mathbf z}}
\newcommand{\bba}{{\mathbf a}}
\newcommand{\bbb}{{\mathbf b}}
\newcommand{\bbc}{{\mathbf c}}


\title{A homotopy method for computing the largest eigenvalue of an irreducible nonnegative tensor}
\author{ Liping Chen\thanks{Internaltional Technology, Amazon Corporate LCC, Seattle, WA 98109, USA. Email: \texttt{lipinc@amazon.com}}  \and
Lixing Han\thanks{Department of Mathematics, University of Michigan-Flint, Flint, MI 48502, USA. Email: \texttt{lxhan@umflint.edu}} 
\and Hongxia Yin \thanks{Department of Mathematics and Statistics, Minnesota State University, Mankato, MN 56001, USA. Email: \texttt{hongxia.yin@mnsu.edu}}
\and Liangmin Zhou\thanks{Mathematics Department, Seattle University, Seattle, WA 98122, USA. Email: \texttt{zhouli@seattleu.edu}}
}
\date{November 6, 2016}
\maketitle

\begin{abstract}  
 In this paper we propose a homotopy method to compute the largest eigenvalue and a corresponding eigenvector of a nonnegative tensor. We prove that it converges to the desired eigenpair when the tensor is irreducible. We also implement the method using an prediction-correction approach for path following. Some numerical results are provided to illustrate the efficiency of the method.
 
\end{abstract}

\ \\
{\bf Key words.}  nonnegative tensors, eigenvalues, Perron pair, homotopy method.

\ \\
{\bf AMS subject classification (2010).} 65F15,  65H17, 65H20.

\section{Introduction}
\label{Intro}
\setcounter{equation}{0}
Let $\bR$ and $\bC$ be the real and complex fields respectively. We denote the set of all $m$th-order, $n$-dimensional real tensors  by $\bR^{[m,n]}$.  
For a tensor $\cA \in \bR^{[m,n]}$ and a vector $x \in \bC^n$, let $\cA x^m$ denote the multilinear form
$$
\cA x^m = \sum_{i_1, \cdots, i_m =1}^{n} A_{i_{1}  \cdots  i_m}x_{i_1} \cdots x_{i_m},
$$
 $\cA x^{m-1}$ denote a column vector in $\bC^n$ whose $i$th entry is
$$
(\cA x^{m-1})_i = \sum_{i_2, \cdots, i_m =1}^{n} A_{i  i_{2}  \cdots  i_m}x_{i_2} \cdots x_{i_m},
$$ 
and $\cA x^{m-2}$ denote an $n \times n$ matrix whose $(i,j)$ entry is
$$
(\cA x^{m-2})_{ij} = \sum_{i_3, \cdots, i_m =1}^{n} A_{i j  i_{3}  \cdots  i_m}x_{i_3} \cdots x_{i_m}.
$$ 

\begin{10} [\cite{CPZ08, Lim05, Qi05}]
{\rm
We say that $(\lambda, x) \in \bC \times (\bC^{n} \backslash \{0\})$ is an eigenpair (eigenvalue-eigenvector) of $\cA \in \bR^{[m,n]}$ if  
\be
\label{evalue}
\cA x^{m-1} = \lambda  x^{[m-1]}, 
\ee
where $x^{[\alpha]} = [x^{\alpha}_{1}, x^{\alpha}_{2}, \ldots, x^{\alpha}_{n}]^T $ for real number $\alpha$. 
}
\end{10}

A tensor $\cA$ is nonnegative (positive) if each of its entries is nonnegative (positive). We denote the set of all $m$th-order, $n$-dimensional nonnegative and positive tensors by $\bR_{+}^{[m,n]}$ and $\bR_{++}^{[m,n]}$ respectively.  The celebrated Perron-Frobenius theory for nonnegative matrices has recently been extended to nonnegative tensors (see for example, \cite{CPZ08, FGH13, YY10, YY11}).  Equipped with this theory, the largest eigenvalue of a nonnegative tensor $\cA$ and its corresponding eigenvector play an important role in various applications, such as in the spectral hypergraph theory, higher order Markov chains, and automatic control, to name a few.  We recall the concept of irreducibility of a tensor:   

\begin{10} [\cite{CPZ08, Lim05}]
{\rm
Denote $N=\{1,2,\ldots,n\}$. A tensor $\cA \in \bR^{[m,n]} $ is said to be  {\it reducible}, if there is a nonempty proper subset $I$ of $N$ such that 
$$
A_{i_1 i_2 \ldots  i_m} =0, \ \ \  \forall i_1 \in I, \   \forall i_2, \ldots, i_m \in N\backslash I.
$$
We say that $\cA$ is {\it irreducible} if it is not reducible. 
}
\end{10}  

We define the spectral radius of a tensor $\cA \in \bR^{[m,n]}$ as 
$$
\rho(\cA) = \max \{ |\lambda|: \lambda \ {\rm is} \ {\rm an} \ {\rm eigenvalue} \ {\rm of} \ \cA\}. 
$$
We also denote $\lambda_*=\rho(\cA)$.  Then we have the following theorem.

\begin{2}  [\cite{CPZ08}]
\label{PF}
If $\cA \in \bR_{+}^{[m,n]}$,   $\lambda_* $ is an eigenvalue of $\cA$ and it has  a   nonnegative corresponding eigenvector $x_*$. \\
If furthermore $\cA$ is irreducible, then  the following hold: 
\begin{itemize}
\item  $\lambda_*>0$. 
\item  Each entry of $x_*$ is positive. 
\item If $\lambda$ is an eigenvalue with a nonnegative eigenvector, then $\lambda = \lambda_*$. Moreover, the nonnegative eigenvector is unique up to a constant multiple. 

\end{itemize}
\end{2}

Under normalization $x_{*}^Tx_{*}=1$, the eigenvector $x_*$ is unique. In this case, we follow \cite{LGL16} and call $\lambda_*$, $x_*$, and $(\lambda_*, x_*)$ the Perron value, the Perron vector, and the Perron pair for the irreducible nonnegative tensor $\cA$, respectively.

Several algorithms have been proposed for finding the Perron pair of a nonnegative tensor in the literature. Ng, Qi, and Zhou \cite{NQZ09} proposed a power-type method. The convergence of this method is established in \cite{CPZ11} for primitive tensors and in \cite{FGH13} for weakly primitive tensors.  Modified versions of the NQZ method have been proposed in \cite{LZI10, ZQX12, ZQW13}.  In \cite{ZQ12}, linear convergence of the NQZ method is proved for essentially positive tensors using a particular starting vector.

To achieve faster convergence, Ni and Qi \cite{NQ15} proposed using Newton's algorithm to solve a polynomial system. Their algorithm is proved to be locally quadratically convergent when the nonnegative tensor is irreducible. They also globalized Newton's method by incorporating a linear search.  However, their globalized algorithm is proved to converge to a stationary point. Thus, there is no guarantee the found solution is the Perron pair.   Liu, Guo, and Lin \cite{LGL16} recently proposed an algorithm that combines Newton's and Noda's iterations for third order nonnegative tensors. This algorithm preserves positivity and is shown to be quadratically convergent to the Perron pair for irreducible nonnegative tensors. It is yet to be seen if this algorithm can be extended to tensors of order $m \geq 4$.

According to Theorem \ref{PF}, finding the Perron pair of an irreducible nonnegative tensor is equivalent to  solving  the system of polynomials (\ref{evalue}) for a particular positive solution $(\lambda_*, x_*)$. An attractive class of methods for solving polynomial systems is the homotopy continuation methods, see for example, \cite{AG90, Li15, SW05}.   Recently, the homotopy techniques have been successfully used to compute generalized eigenpairs of a tensor in \cite{CHZ15, CHZ16}). However, the homotopy methods of \cite{CHZ15, CHZ16} are not suitable for finding the Perron pair of a large size irreducible  nonnegative tensor because those methods are designed to find all (real and complex) eigenpairs.  
In this paper, we propose a homotopy method for computing the Perron pair of an irreducible nonnegative tensor.   This method is suitable for large size tensors.  The iterations in our method follow a curve in the positive orthant  $\bR_{++}^{n+1}$. Thus it preserves positivity.

The paper is organized as follows.  We  introduce our homotopy method and prove its convergence in Section 2. Then we describe an implementation of the method and give some numerical results in Section 3.  Some concluding remarks are given in Section 4.

\section{A Homotopy Method}
\label{sec2}
\setcounter{equation}{0}

Our goal is to  compute the Perron pair  of an irreducible nonnegative tensor $\cA$ using a homotopy method.  For this purpose, we will solve the following problem:
\be
\label{target}
P(\lambda,x) = 
\begin{pmatrix}
\cA x^{m-1} - \lambda x^{[m-1]}\\
x^Tx -1
\end{pmatrix} = 0.
\ee
Choose  positive vectors $a$ and $b$ from $  \bR_{++}^{n}$. Define the positive tensor $\cE \in \bR^{[m,n]}_{++}$ such that 
\be
\label{starttensor}
\cE= a^{[m-1]} \circ b \circ \cdots \circ b \in \bR^{[m,n]}_{++}
\ee
where $\circ$ denotes the outer product. Clearly the 
$(i_1, i_2,\ldots, i_m)$ entry of $\cE$ is
$$
E_{i_1i_2 \ldots i_m}= a_{i_1}^{m-1} b_{i_2} \cdots b_{i_m}. 
$$
For the tensor $\cE$, we have 
\begin{1}
\label{lemma1}
For any $a, b \in \bR_{++}^{n}$, the tensor $\cE$ is a positive tensor in $\bR^{[m,n]}_{++}$. Moreover, the Perron pair  $(\lambda_0, x_0)$ of $\cE$ satisfies
\be
\label{beginpair}
\lambda_0=(a^Tb)^{m-1},  \ x_0=a/\|a\|.
\ee
\end{1}

Using the start system
\be
\label{start}
Q(\lambda,x) = 
\begin{pmatrix}
\cE x^{m-1} - \lambda x^{[m-1]}\\
x^Tx -1
\end{pmatrix} = 0,
\ee
we construct the following homotopy 
\be
\label{homotopy}
H(\lambda,x,t) = (1-t) Q(\lambda,x) + t P(\lambda,x) =0,  \ \ \   t \in [0,1].
\ee
Clearly we have
\be
\label{hexpress}
H(\lambda,x,t) = 
\begin{pmatrix}
(t \cA + (1-t)\cE) x^{m-1} - \lambda x^{[m-1]}\\
x^Tx -1
\end{pmatrix}.
\ee

The Jacobian matrix of $H(\lambda, x, t)$ plays an important role in designing a homotopy method for solving the eigenvalue problem (\ref{evalue}). In order to compute this matrix, we need to partially symmetrize  tensor $\cA=(A_{i_1,i_2,\ldots, i_m})$ for indices $i_2, \ldots, i_m$. Specifically, we define the partially symmetrized tensor 
$\bar{\cA}=(\bar{A}_{i_1, i_2, \ldots, i_m})$ by
\be
\label{psymm}
\bar{A}_{i_1 i_2 \ldots i_m} = \frac{1}{(m-1)!} \sum_{\pi} A_{i_1 \pi(i_2\ldots   i_m)},
\ee
 where the sum is over all the permutations $\pi(i_2\ldots  i_m)$. The Jacobian matrix of $\cA x^{m-1}$ with respect to $x$ is 
\be
\label{jacobian}
D_x \cA x^{m-1} = (m-1) \bar{\cA} x^{m-2}. 
\ee
Note that the tensor $\cE$ is partially symmetric. Thus the partially symmetrized tensor $\bar{\cB}_t$ of tensor $\cB_t = t \cA + (1-t)\cE$ is 
given by 
$$
\bar{\cB}_t = t \bar{\cA} + (1-t)\cE. 
$$

The partial derivatives of $H$ with respect to $\lambda$, $x$, and $t$ are:
$$
D_{\lambda} H (\lambda,x,t)  = \begin{pmatrix}
-  x^{[m-1]}\\
 0
\end{pmatrix},
$$
$$
D_{x} H (\lambda,x,t)  = \begin{pmatrix}
(m-1) [ \bar{\cB_t} x^{m-2} - \lambda C]     \\
  2x^T
\end{pmatrix},
$$
and 
$$
D_{t} H (\lambda,x,t)  = \begin{pmatrix}
 (\cA-\cE)x^{m-1}\\
 0
\end{pmatrix},
$$
respectively, where   $C$ is the diagonal matrix
$$
C = {\rm diag} ([x_{1}^{m-2}, x_{2}^{m-2}, \ldots, x_{n}^{m-2}])
$$

\begin{1}
\label{lemma2} Suppose that $\cA \in \bR_{+}^{[m,n]}$  and $a, b \in \bR_{++}^n$.  Then we have the following: \\
(a) For any $t \in [0,1)$, $H(\lambda, x, t)=0$ has a unique solution $(\lambda(t), x(t))$ in $\bR_{++} \times \bR_{++}^{n}$, which is 
the Perron pair of the positive tensor $t\cA+(1-t)\cE$. Moreover,   
 the Jacobian matrix 
$$
D_{(\lambda,x)} H(\lambda(t), x(t), t) = 
 \begin{pmatrix}
-  x(t)^{[m-1]}   & (m-1) [ \bar{\cB_t} x(t)^{m-2} - \lambda(t) C(t)]     \\
 0 &  2x(t)^T
\end{pmatrix}
$$ 
is of rank $n+1$, where $C(t)= {\rm diag} ([x_{1}(t)^{m-2}, x_{2}(t)^{m-2}, \ldots, x_{n}(t)^{m-2}])$.  \\
(b) If furthermore $\cA$ is irreducible, then $H(\lambda, x, 1)=P(\lambda, x)=0$ has a unique solution  $(\lambda_*,x_*)=(\lambda(1), x(1))$ in $\bR_{++} \times \bR_{++}^{n}$, which is 
the Perron pair of the  tensor $\cA$.
Moreover, the Jacobian matrix 
$$
D_{(\lambda,x)} H(\lambda_*, x_*, 1) = 
 \begin{pmatrix}
-  x_*^{[m-1]}   & (m-1) [ \bar{\cA} x_*^{m-2} - \lambda_* C_*]     \\
 0 &  2x_{*}^T
\end{pmatrix}
$$ 
 is of rank $n+1$, $C_*= {\rm diag} ([x_{*1}^{m-2}, x_{*2}^{m-2}, \ldots, x_{*n}^{m-2}])$.
\end{1}

\bp
For part (a), by Theorem \ref{PF}, $(t\cA+(1-t)\cE) x^{m-1} - \lambda x^{[m-1]} =0 $ has a unique solution in  $\bR_{++} \times \bR_{++}^{n}$, up to a constant multiple of $x$. 
By imposing the normalization condition $x^Tx=1$, $P(\lambda, x)=0$ has a unique solution $(\lambda(t), x(t))$ in $\bR_{++} \times \bR_{++}^{n}$. Clearly  $(\lambda(t), x(t))$ is the Perron pair of  tensor $t\cA+(1-t)\cE$. Now by \cite[Lemma 4.1]{NQ15}, the Jacobian matrix $ D_{(\lambda,x)} H(\lambda(t), x(t), t)$ is nonsingular, that is, it has rank $n+1$.  \\
The proof of part (b) is similar. 
\ep

\begin{1}
\label{lemma3}
Suppose that $\cA \in \bR_{+}^{[m,n]}$  and $a, b \in \bR_{++}^n$. Then 
 $H_{+}^{-1}(0) = \{ (\lambda, x, t) \in  \bR_{++} \times \bR_{++}^{n} \times [0,1): H(\lambda, x, t)=0 \} $ is  a one-dimensional smooth manifold. \\
\end{1}
\bp
The conclusion follows from Lemma \ref{lemma2} and the implicit function theorem. 
\ep

\begin{1}
\label{lemma4}
Suppose that $\cA \in \bR_{+}^{[m,n]}$  and $a, b \in \bR_{++}^n$. The $H_{+}^{-1}(0)$  in Lemma \ref{lemma3} is uniformly bounded for $t \in [0,1)$. 
\end{1}
\bp
By the Gerschgorin  theorem for the eigenvalues of tensors (\cite{Qi05, CQZ13}), for the eigenvalue $\lambda$ of tensor $t\cA+(1-t)\cE$,
$$
|\lambda| \leq \sum_{i_1=1, \ldots, i_m=1}^{n} |(t\cA+(1-t)\cE)_{i_1i_2 \ldots i_m}| \leq \sum_{i_1=1, \ldots, i_m=1}^{n} |A_{i_1i_2 \ldots i_m}| + \sum_{i_1=1, \ldots, i_m=1}^{n} |E_{i_1i_2 \ldots i_m}|.
$$
Note that $\|x\|=1$. Thus, we have
$$
\| (\lambda,x,t) \| \leq |\lambda|+\|x\|+|t| \leq 2+ \sum_{i_1=1, \ldots, i_m=1}^{n} |A_{i_1i_2 \ldots i_m}| + \sum_{i_1=1, \ldots, i_m=1}^{n} |E_{i_1i_2 \ldots i_m}|.
$$
\ep

Our main theorem shows that the homotopy (\ref{homotopy}) works.

\begin{2}
\label{main}
Suppose that $\cA \in \bR_{+}^{[m,n]}$  and $a, b \in \bR_{++}^n$. Starting from the Perron pair $(\lambda_0,x_0)$ of $\cE$ as defined in (\ref{beginpair}),  let $(\lambda(t), x(t))$ be the solution curve obtained by solving the homotopy $H(\lambda,x,t)=0$ in $\bR_{++} \times \bR_{++}^{n} \times [0,1)$. Then each limit point of $(\lambda(t), x(t))$ is a nonnegative eigenpair of $\cA$. \\
 If furthermore $\cA$ is irreducible, then 
\be
\label{convergence}
    \lim_{t \to 1} (\lambda(t), x(t)) = (\lambda_*, x_*),
\ee
where $(\lambda_*, x_*)$ is the Perron pair of $\cA$. 
\end{2}
\bp
According to Lemma \ref{lemma3}, $H_{+}^{-1}(0) = \{ (\lambda, x, t) \in \bR_{++} \times \bR_{++}^{n} \times [0,1): H(\lambda, x, t)=0 \} $ is a one-dimensional smooth manifold  in $\bR_{++} \times \bR_{++}^{n}$.  Since $H_{+}^{-1}(0)$ is uniformly bounded by Lemma \ref{lemma4}, $(\lambda(t), x(t))$ is uniformly bounded on $[0,1)$ and thus its limit set is not empty when $t \to 1$. Let $(\lambda_*, x_*)$ be a limit point of $(\lambda(t), x(t))$ as $t \to 1$.  Then $(\lambda_*, x_*) \in \bR_{+} \times \bR_{+}^{n} $ and 
\be
\label{pfeqn}
\begin{pmatrix}
 \cA x_{*}^{m-1} - \lambda x_{*}^{[m-1]} \\
x_{*}^Tx_* -1
\end{pmatrix} = 0,
\ee
by taking the limit in (\ref{hexpress}). This $(\lambda_*,x_*) $ is a nonnegative eigenpair of $\cA$. \\
 If in addition  $\cA$ is irreducible, then (\ref{pfeqn}) implies that  $(\lambda_*, x_*)$ must be the unique Perron pair of $\cA$. This further implies that (\ref{convergence}) holds.
\ep

\begin{5}{\rm
To follow the curve in the homotopy method, we differentiate $H(\lambda,x, t)=0$ with respect to  $t$. Then
\be
\label{diffeq}
D_{(\lambda,x)}  H    \cdot \left ( \begin{array}{c}
 \frac{d \lambda }{d t} \\
\frac{d x}{d t} 
\end{array} 
\right ) 
= - D_t H.
\ee
Since $D_{(\lambda,x)}  H $ is nonsingular for all $t \in [0,1)$, this system of differential equations is well defined. 
We will follow the curve by solving this system with the initial condition $(\lambda(0), x(0))=(\lambda_0,x_0)$, where $\lambda_0, x_0$ are defined in (\ref{beginpair}). 
}
\end{5}

\section{Implementation and Numerical Results}
\label{sec3}
\setcounter{equation}{0}

We now present an algorithm that implements the homotopy method proposed in the previous section for computing the Perron pair of an irreducible nonnegative tensor $\cA \in \bR^{[m,n]}_{+}$.  An Euler-Newton type predication-correction approach for solving the system of differential equations (\ref{diffeq}) with  initial conditions $(\lambda(0),x(0))=(\lambda_0,x_0)$ 
is  used.

\begin{70}
\label{algorithmA}

\begin{framed}

\ \\
{\bf  \large Initialization.}  Choose positive vectors $a, b \in \bR_{++}^{n}$ and construct the tensor $\cE= a^{[m-1]} \circ b \circ \cdots \circ b  \in \bR^{[m,n]}_{++}$. Choose initial step size $\Delta t_0>0$, tolerances $\epsilon_1>0$ and $\epsilon_2>0$. Let 
$t_0=0$, $x_0=a/\|a\|$, $\lambda_0=(a^Tb)^{m-1}$. Set $k=0$. \\
\ \\
{\bf  \large Path following.} 
For $k=0,1,\cdots$ until $t_N<1$ and $t_{N+1} \geq 1$ for some $N$: \\
  Set $ t_{k+1} = t_k+ \Delta t_k$. If $t_k<1$ and $t_{k+1} \geq 1$, then set $N=k$ and reset $t_{N+1} = 1$ and $\Delta t_N=1-t_N$.   \\
 Let $u=(\lambda,x)$ and $u_k= (\lambda_k, x_k)$. To find the next point on the path $H(u,t) = 0$, we employ the following Euler-Newton prediction-correction strategy:
\begin{itemize}
\item Prediction Step: Compute the tangent vector $\dfrac{du}{dt}$ to $H(u,t)=0$ at $t_k$ by solving the linear system
\begin{equation*}
D_u H(u_k,t_k)\dfrac{du}{dt} = -D_tH(u_k,t_k)
\end{equation*}
for $\dfrac{du}{dt}$. Then compute the approximation 
$\tilde{u}$ to $u_{k+1}$ by
$$
\tilde{u} = u_k + \Delta t_k \frac{du}{dt}.
$$
\item Correction Step: Use Newton's iterations. Initialize $v_0 = \tilde{u}$. For $i = 0, 1, 2, \dots$, compute
$$
v_{i+1} = v_i - [D_u H(v_i,t_{k+1})]^{-1} H(v_i,t_{k+1})
$$
until $\|H(v_{J},t_{k+1})\| \leq \epsilon_1$ if $k<N$ or $\|H(v_{J},t_{k+1})\| \leq \epsilon_2$ if $k=N$. Then let $u_{k+1} = v_{J}$. If $k=N$, we set $(\lambda_*, x_*)= u_{N+1} $ as the computed pair and stop.
 \item Adaptively updating the step size $\Delta t_k$:  If more than three steps of Newton iterations were required to converge within the desired accuracy, then $\Delta t_{k+1}=0.5 \Delta t_k$. If $\Delta t_{k+1} \leq 10^{-6}$, set $\Delta t_{k+1}=10^{-6}$.  If two consecutive steps were not cut, then  $\Delta t_{k+1}=2 \Delta t_k$. If $\Delta t_{k+1} \geq 0.5$, set $\Delta t_{k+1}=0.5$. Otherwise, $\Delta t_{k+1}=\Delta t_k$.  
\end{itemize}
 Set $k=k+1$.  
\end{framed}
\end{70}

We now report some numerical results to test Algorithm \ref{algorithmA}. We compare it with the NQZ method \cite{NQZ09}. All the experiments were done using MATLAB 2014b on a laptop computer with Intel Core i7-4600U at 2.10 GHz and 8 GB memory running Microsoft Windows 7.  The tensor toolbox of \cite{BK15} was used to compute tensor-vector products and to compute partially symmetrized tensor $\bar{\cA}$.  

In our experiments, we used $a=b=[1,1,\ldots,1]^T$ and $\Delta t_0=0.1$, $\epsilon_1=10^{-5}$ and $\epsilon_2=10^{-12}$ in Algorithm \ref{algorithmA} and the initial vector $x_0 = [1,1,\ldots,1]^T/\sqrt{n}$ in the NQZ method.  
We also preprocessed the tensor $\cA$ by setting 
\be
\label{preprocess}
\tilde{\cA} = \cA / \tau,
\ee  
where $\tau:=\max(\cA)$ is the largest entry in $\cA$. After finding the Perron pair $(\tilde{\lambda}_*, x_*)$ by using Algorithm \ref{algorithmA} or the NQZ method on $\tilde{\cA}$, the Perron pair for $\cA$ is obtained by setting $(\lambda_*,x_*)=(\tau \tilde{\lambda}_*,x_*)$. The NQZ method was terminated if one of the following conditions was met:

(a). $\|[\tilde{\cA} x_{k}^{m-1} - \tilde{\lambda_k} x_{k}^{[m-1]} ; x_{k}^Tx_k - 1]\| \leq 10^{-12}$.

(b). The number of iterations exceeds 10000. 

Note that regular termination condition (a) is the same as the one used in Algorithm \ref{algorithmA} at regular termination. We also set the maximal allowed number  of prediction-correction steps for Algorithm \ref{algorithmA} as 10000.   

Our first example is adapted from \cite[Example 3.6]{CPZ11}. 

\begin{60}
\label{example1} {\rm
 Consider  $\cA \in \bR^{[3,3]}_{+}$ defined by:
$$
a_{122}=1, a_{133}=2, a_{211}=3, a_{311}=4, \ {\rm and} \ a_{ijk}=0 \ {\rm otherwise}.
$$
}
\end{60}

This tensor is irreducible but not primitive. Algorithm \ref{algorithmA} converged to the Perron pair using 5 prediction-correction steps and 14 Newton iterations, while the NQZ method did not converge after 10000 iterations.  We remark that the NQZ method converges to the Perron pair if the strategy of adding a shift to $\cA$ introduced in \cite{LZI10}  is used. For example,  applying the NQZ method to the tensor $\cA  + \cI$, the NQZ method can find the Perron pair in 25 iterations.

\begin{60}
\label{example2}
{\rm
 Consider $\cP \in \bR^{[3,3]}_{++}$ from \cite[Example 1]{LGL16}, which is defined by
\[
\cP(1,:,:) = \left [
\begin{array}{ccc}
0.9000 & 0.6700 & 0.6604 \\
0.3340 & 0.1040 & 0.0945 \\
0.3106 & 0.0805 & 0.0710 
\end{array} 
\right ],
\]
\[
\cP(2,:,:) = \left [
\begin{array}{ccc}
0.0690 & 0.2892 & 0.0716 \\
0.6108 & 0.8310 & 0.6133 \\
0.0754 & 0.2956 & 0.0780 
\end{array} 
\right ],
\]
\[
\cP(3,:,:) = \left [
\begin{array}{ccc}
0.0310 & 0.0408 & 0.2680 \\
0.0552 & 0.0650 & 0.2922 \\
0.6140 & 0.6239 & 0.8510 
\end{array} 
\right ].
\]
We set tensor $\cA = \cP + \gamma \cI$, where $ \cI \in \bR^{[3,3]} $ is the identity tensor and $\gamma$ is a parameter.  We summarize the numerical results for Example \ref{example2} in  Table \ref{table1}. In this table and Table \ref{table2}, for Algorithm A, \texttt{itr} and \texttt{nwtitr} denote the number of prediction-correction steps and the total number of Newton iterations were used, respectively. For the NQZ method, \texttt{itr} denotes the number of iterations was used. \texttt{time} represents the CPU time used  when the termination condition (a) was met.  
}
\end{60}

\begin{table}[htbp]
\begin{center}
\begin{tabular}{|c|c|c|c|| c|c|c|}
\hline
& \multicolumn{3}{|c||}{Algorithm \ref{algorithmA}} & \multicolumn{2}{|c|}{NQZ} \\ \hline
$\gamma$ & itr & nwtitr & time & itr & time \\ \hline
0 & 5 & 15 & 0.1255 	&   29 & 0.0345 \\ \hline
10 & 5 & 13 & 0.1237 	&   157 & 0.1878 \\ \hline
$10^2$ & 5 & 13 & 0.1248 	 & 1195 & 1.2473 \\ \hline
$10^3$ & 5 & 11 & 0.1019 	&  $>10000$  &  \\ \hline
$10^4$ & 5 & 9 & 0.1003 	&   $>10000$ &  \\ \hline
\end{tabular}
\end{center}
\caption{Performance of Algorithm \ref{algorithmA} and the NQZ method on Example \ref{example2}}
\label{table1}
\end{table}

\begin{60}
\label{example3}
{\rm
 Let  $\cC \in \bR^{[m,n]}_{+}$ be a tensor, each of  its entries being a random number uniformly distributed in $[0,1]$.  Define
$$
\cA = \cC + \gamma \cI,
$$
where $\cI$ is the identity tensor in $\bR^{[m,n]}$ and $\gamma$ is a parameter.  The numerical results for this example are given in   Table \ref{table2}.
}
\end{60}

\begin{table}[htbp]
\begin{center}
\begin{tabular}{|c|c|c|c|c|| c|c|c|}
\hline
&& \multicolumn{3}{|c||}{Algorithm \ref{algorithmA}} & \multicolumn{2}{|c|}{NQZ} \\ \hline
$(m,n)$ &$\gamma$ & itr & nwtitr & time  & itr & time \\ \hline
(3,20) & $10^2$  & 5 & 12 & 0.1352  & 22  & 0.0519  \\ \hline
(3,20) &  $10^4$  &5  & 10  & 0.1261  & 942  & 1.0130  \\ \hline
(3,20) & $10^6$  &  5 & 8  & 0.1167 & $>10000$  &  \\ \hline
(3,100) & $10^2$  & 5  & 12  & 0.3258  & 7  & 0.0506  \\ \hline
(3,100) & $10^4$ & 5 & 12  & 0.3265  & 49  & 0.1318 \\ \hline
(3,100) & $10^6$ & 5  & 8 & 0.2924 & 2998  & 5.7511  \\ \hline
(3,200) & $10^2$ & 5  & 12 & 1.1023  & 7 & 0.0915  \\ \hline
(3,200) &$10^4$  & 5 & 12  & 1.0684 & 19  & 0.1643  \\ \hline
(3,200) & $10^6$ & 5 & 8 & 0.9900  & 774 &  4.1408 \\ \hline
(3,200) & $10^7$ & 5  & 8 & 0.9931  & 6510  & 35.8417  \\ \hline
(4,10) & $10^2$ & 5  & 12 & 0.1331 &13  & 0.0505  \\ \hline
(4,10) & $10^4$ & 5  & 10 & 0.1241  & 361  & 0.4277  \\ \hline
(4,10) & $10^6$ & 5  & 8 & 0.1240  & $>10000$  &  \\ \hline
(4,50) & $10^2$ & 5 & 12 & 1.1523 & 4 & 0.0842  \\ \hline
(4,50) & $10^4$ &  5 & 9 & 1.1133  & 10  & 0.1110  \\ \hline
(4,50) & $10^6$  & 5 & 7  & 1.0525 & 244 & 1.3374  \\ \hline
(4,50) & $10^7$  & 5 & 7 & 1.0465 & 2004  & 10.9171  \\ \hline
(4,100) & $10^2$ & 5  & 12  & 22.6062  & 4 & 0.8868  \\ \hline
(4,100) & $10^4$ & 5 & 11 & 21.2236 & 5  & 0.8900  \\ \hline
(4,100) & $10^6$ & 5 & 7 & 20.2311 & 37  & 3.4189 \\ \hline
(4,100) & $10^8$ & 5 & 7 & 20.4364  & 2078 & 168.9237 \\ \hline

\end{tabular}
\end{center}
\caption{Performance of Algorithm \ref{algorithmA} and the NQZ method on Example \ref{example3}}
\label{table2}
\end{table}

From Tables \ref{table1} and \ref{table2}, we observe that Algorithm \ref{algorithmA} is more efficient than the NQZ method when the parameter $\gamma$ is large, in terms of numbers of iterations. This is because the NQZ is linearly convergent and its  rate of convergence depends on the ratio $r/\lambda_*$, where $r$ is the maximum modulus of the eigenvalues of $\cA$ distinct from $\lambda_*$. When $\gamma$ is large, the ratio  $r/\lambda_*$ is close to $1$ (For example, this ratio is $0.9988$ for Example \ref{example2} when $\gamma=10^3$). Thus the NQZ method becomes slow. Using the shifting strategy introduced in \cite{LZI10} does not improve the performance of the NQZ method for these examples when $\gamma$ is large. 

On the other hand, the last correction step in Algorithm \ref{algorithmA} is Newton's method. Therefore, it is quadratically convergent.  Numerical results show that the ratio $r/\lambda_*$ does not affect the performance of Algorithm \ref{algorithmA}. We remark that  for the $m=4, n=100$ cases, the relatively large CPU time used by Algorithm \ref{algorithmA} is mainly due to the procedure of partially symmetrizing tensor $\cA$.  A more efficient symmetrization method can help save the CPU time for such cases. Of course,  the partial symmetrization is not needed if the tensor $\cA$ is symmetric.

\section{Concluding Remarks}
\label{sec4}

We have proposed a homotopy method for computing the largest eigenvalue and a corresponding eigenvector of a nonnegative tensor. We have proved that this method converges to the Perron pair for irreducible nonnegative tensors. We have implemented it using an Euler-Newton prediction-correction approach for path following in Algorithm \ref{algorithmA}.  Our numerical results show the efficiency of Algorithm \ref{algorithmA}. This algorithm is promising when the ratio $r/\lambda_*$ is close to $1$, where $r$ the maximum modulus of the eigenvalues distinct from the Perron value $\lambda_*$. In this situation, the Power-type methods can be slow.

The algorithms proposed in this paper and in \cite{CHZ15,CHZ16} show that the homotopy techniques are very useful for computing tensor eigenvalues.  An interesting direction for  future research  is to apply a homotopy approach to compute extreme Z-eigenvalues for nonnegative tensors.  Another direction is to  employ this approach to compute extreme eigenvalues of symmetric tensors.

\ \\
{\bf Acknowledgements.}  Lixing Han was supported in part by a 60th Anniversary Research Grant, Office of the Provost, UM-Flint.


\end{document}